\documentclass[11pt]{elsarticle}
\usepackage{amssymb}
\usepackage{amsmath}
\usepackage{amsfonts,graphicx,color, bbm}
\usepackage{mathtools}
\newtheorem{definition}{Definition}[section]
\newtheorem{theorem}{Theorem}[section]
\newtheorem{lemma}{Lemma}[section]
\newtheorem{corollary}{Corollary}[section]
\newtheorem{proposition}{Proposition}[section]
\newtheorem{example}{Example}[section]
\newtheorem{remark}{Remark}[section]
\newcommand{\ignore}[1]{}{}
\newproof{pf}{Proof}
\newproof{pot2}{Proof of Theorem \ref{Thm2}}
\newproof{pot3}{Proof of Proposition \ref{prop-yinna1}}

\newcommand\beq{\begin{equation}}
	\newcommand\eeq{\end{equation}}

\makeatletter
\newcommand{\rmnum}[1]{\romannumeral #1}
\newcommand{\Rmnum}[1]{\expandafter\@slowromancap\romannumeral #1@}

\def\geq{\geqslant}
\def\leq{\leqslant}
\def\P{\mathbb{P}}
\def\E{\mathbb{E}}
\def\N{\mathbb{N}}
\def\R{\mathbb{R}}

\newtheorem{hyp}{Hypothesis}
 
\begin{document}
	\let\today\relax
\date{\quad}
\begin{frontmatter}



\title{Convergence Rate of Euler--Maruyama Scheme to the Invariant Probability Measure Under Total Variation Distance for the SDEs}
\author[cor1]{Yuke Wang}
\author[cor2]{Yinna Ye\corref{c1}}

\cortext[c1]{Corresponding author}
\address[cor1]{{School of Mathematics and Statistics, Northeastern University at Qinhuangdao, Qinhuangdao 066004, China}}	
\address[cor2]{Department of Applied Mathematics, School of Mathematics and Physics, Xi'an Jiaotong-Liverpool University, Suzhou 215123, China}	


\begin{abstract}
This article shows the geometric decay rate of the Euler--Maruyama scheme for a one-dimensional stochastic differential equation towards its invariant probability measure under total variation distance. Firstly, the existence and uniqueness of invariant probability measure and the uniform geometric ergodicity of the chain are studied through the introduction of non-atomic Markov chains. Secondly, the equivalent conditions for uniform geometric ergodicity of the chain are discovered by constructing a split Markov chain based on the original Euler--Maruyama scheme.
\end{abstract}




\end{frontmatter}
\date{\today}

\section{Introduction}
Consider the following stochastic differential equation (SDE) on $\R$:
\begin{equation}\label{SDE}
	dX_t = g(X_t)\,dt + \sigma\, dB_t,\quad X_0 = x_0,
\end{equation}
where $(B_t)_{t\geq0}$ is a standard Brownian motion, $\sigma>0$ is a constant, and~$g: \R\rightarrow\R$
satisfies the hypothesis below. 
\begin{hyp}\label{hyp1}
	There exist $L$, $K_1 > 0$ and $K_2 \geq 0$ such that for every $(x,y)\in\R^2$,
	\begin{align}
		|g(x)-g(y)| &\leq L|x-y|,\label{assump1-1}\\
		(g(x)-g(y))(x-y)&\leq-K_1(x-y)^2+K_2.\label{assump1-2}
	\end{align}
	
	Moreover, $g$ is second-order differentiable, and~the second-order derivative of $g$ is bounded.
\end{hyp}	

The inequality (\ref{assump1-1}) implies for any $x\in\R$,
\begin{equation}\label{eq-g-square}
	g^2(x)\leq 2L^2x^2+2{g^2(0)},\quad |g'(x)|\leq L,
\end{equation}
where $g'$ represents the derivative of the function $g$.
Moreover, the~inequality (\ref{assump1-2}) and Young's inequality imply
\begin{equation}\label{assump1-3}
	xg(x)\leq -\frac{K_1}{2}x^2+c_g,
\end{equation}
with $c_g=K_2+\frac{1}{2K_1}g^2(0)$.

For a given step size $\eta\in(0,1)$, the~Euler--Maruyama (EM) scheme $(\theta_k)_{k\geq0}$ for the SDE  (\ref{SDE}) is given by the following recursive relation: for any $k\geq0$,
\begin{equation}\label{recursive}
	\theta_{k+1}=\theta_{k}+\eta\, g(\theta_k)+\sqrt{\eta}\,\sigma\varepsilon_{k+1},
\end{equation}
where $(\varepsilon_k)_{k\geq1}$ are independent and identically distributed standard normal random variables.

Over the past decades, EM schemes have been widely and intensively used for solution approximations of SDEs on $\R^d$ with $d\geq1$ (see for instance \cite{BJ22,BPP22,HL18,KS19,LZ24,PT17,P08,PT97}). In~the case when $g(x)=-\nabla U(x)$ with $U$ being a potential, (\ref{SDE}) is called Langevin diffusion, and~(\ref{recursive}) is called unadjusted Langevin algorithm (ULA) with constant step size. The~ULA is nowadays a popular class of Markov chain Monte Carlo (MCMC) algorithms (\cite{BGJM11}), which is commonly applied to solve Monte Carlo sampling problems especially in the field of machine learning (see, for~instance \cite{MY11,DL21}). Roberts and Tweedie~\cite{RT96} initially provided probabilistic analysis of asymptotic behaviour of Langevin diffusion and ULA. They found necessary (Theorem 2.2 {\cite{RT96}}) and sufficient conditions (Theorem 2.3,~\cite{RT96}) for exponential ergodicity of Langevin diffusion, and~sufficient conditions for geometric ergodicity of ULA (Theorem 3.1,~\cite{RT96}) respectively. Dalalyan~\cite{D17} found an upper bound for the error in total variation (TV) distance between $\theta_k$ and the ergodic measure of $(X_t)_{t\geq0}$. Durmus and Moulines~\cite{DM17} established geometric ergodicity in TV distance for Langevin diffusion under similar convexity and Lipschitz conditions, with~explicit convergence rates. Recently, Fan~et~al.~\cite{FHX24} established normalized and self-normalized Cram\'{e}r-type moderate deviations for the EM scheme for the SDE (\ref{SDE}) on $\R^d$ with $d\geq1$. Meanwhile, the~geometric ergodicity of the ULA towards its invariant probability measure under appropriate distance has been rarely studied so far in the literature. As~an extension of classical Markov chains with discrete state space, the~EM scheme studied in this paper is a non-atomic Markov chain, which is also of independent interest for the study on other stochastic models with  a non-independence~structure. 

The objectives of this work are to prove the geometric ergodicity of the EM scheme given by (\ref{recursive}) for the SDE (\ref{SDE}) under TV distance, and~provide furthermore a convergence rate to its invariant probability~measure. 

Let $\mathcal{B}(\R)$ be the Borel $\sigma$-algebra on $\R$. From~(\ref{recursive}), it can be seen that $(\theta_k)_{k\geq0}$ is a Markov chain on the continuous state space $\R$, with~Markov kernel $P_{\eta}$ given by
\begin{align}
	P_{\eta}(x,A)=\frac{1}{\sqrt{2\pi\eta}\,\sigma}\int_A\exp\left(-\frac{(y-x-\eta g(x))^2}{2\eta\sigma^2}\right) dy, \label{density}
\end{align}
for any $x\in\R$ and $A\in \mathcal{B}(\R)$. From~(\ref{density}), $P_{\eta}$ admits a kernel density with respect to (w.r.t) the Lebesgue measure given by $p_{\eta}(x,y):=\frac{1}{\sqrt{2\pi\eta}\,\sigma}\exp\left(-\frac{(y-x-\eta g(x))^2}{2\eta\sigma^2}\right)$. The~Markov kernel $P_\eta$ can be defined equivalently as, for~any $x\in \R$ and measurable function $h$ on $\left(\R,\mathcal{B}(\R)\right)$,
\begin{align*}
	P_{\eta}\,h(x):=\int_{\R}h(y) \,p_{\eta}(x,y)\,dy.
\end{align*}

 Let us introduce now TV norm $\|\cdot\|_{TV}$ and TV distance $d_{TV}$ respectively. Let $(X, \mathfrak{X})$ be a measurable space. Consider $\mathbbm{F}_b$ the set of all bounded and measurable functions from $(X, \mathfrak{X})$ to $\R$. For~any finite signed measure $\xi$ on $(X, \mathfrak{X})$, the~total variation norm of $\xi$ is defined~by
 
 $$\|\xi\|_{TV}:=\sup\left\{\xi(f)\,|\,f\in\mathbbm{F}_b,\,|f|_{\infty}\leq1\right\},$$
 where $|f|_{\infty}:=\sup\left\{|f(x)|\,|\,x\in X\right\}$. For~any two probability measures $\xi$ and $\xi'$ on $(X, \mathfrak{X})$, the~total variation distance between $\xi$ and $\xi'$ is defined by
 \begin{equation}\label{def_TVD}
 	d_{TV}(\xi,\xi'):=\frac{1}{2}\|\xi-\xi'\|_{TV}.
 \end{equation}

Throughout the paper, sometimes supplied with some indices, $\tau$ and $c$ denote positive constants, whose value may vary from line to line. Subsequent symbols with the prefix $\tau$ will be redefined according to specific scenarios and are independent of each other. For~any $(a,b)\in\R^2$, $a\vee b=\max\{a,b\}$ and $a\wedge b=\min\{a,b\}$. When we write ``$C$ is an ($m,\varepsilon\nu$)-small set'', it will be always assumed that $\varepsilon\in(0,1]$ and $\nu$ is a probability~measure.

The paper is organized as follows. Section~\ref{sec2} derives some basic properties (Proposition~\ref{aperiodic}) and existence and uniqueness of the invariant probability measure $\pi_{\eta}$ (Theorem \ref{Thm1}) of the EM scheme $(\theta_k)_{k\geq0}$. Section~\ref{sec3} obtains the convergence rate of the chain $(\theta_k)_{k\geq0}$ towards $\pi_{\eta}$ under TV distance (Theorem \ref{Thm2}), through demonstrating furthermore the equivalent conditions for its uniform geometric ergodicity (Proposition \ref{uniformGE}).  To~prove Proposition  \ref{uniformGE}, we use the splitting construction approach, which is described in Section~\ref{sec4}. In~Section~\ref{sec5}, we discuss some connections and special highlights of the present work with classical ergodic theory for Markov chains on general state space, and~some future research directions. In~the Appendix, some known and useful properties of non-atomic and atomic Markov kernels are collected respectively in
\ref{sec5-1} and \ref{sec5-2}. In \ref{sec5-3}, we state some auxiliary~results.
\section{Existence and Uniqueness of Invariant Probability Measure $ \pi_ \eta$}\label{sec2}
Consider the EM scheme $(\theta_n)_{n\geq0}$ defined by (\ref{recursive}) with Markov kernel $P_{\eta}$ given by $(\ref{density})$. In~the sequel, let $\P_x$ denote the probability measure on the canonical space $\left(\R^{\N},\left(\mathcal{B}(\R)\right)^{\otimes\N}\right)$ of the chain $(\theta_n)_{n\geq0}$, with~the initial state $\theta_0=x$. And~$\E_x$ denotes the expectation under the probability measure $\P_x$. 

\begin{hyp}\label{hyp2}\label{assump2}
	For any $x\in\R$ and $\eta\in(0,1)$, $|x+\eta g(x)|\leq c.$
\end{hyp}

We have the following basic property about the Markov chain $(\theta_n)_{n\geq0}$ and its kernel $P_{\eta}$. Proposition \ref{aperiodic} (1) below provides a sufficient condition for the existence of accessible $(1,\mu)$-small set $C$ satisfying $\mu(C)>0$.
\begin{proposition}\label{aperiodic}
	\begin{enumerate}
		\item[(1)] If $C$ is a compact subset of $\mathbb{R}$ such that $\mbox{Leb}(C)>0$,  then $C$ is an accessible $(1,\epsilon\nu)$-small set, where
		\begin{align}\label{Def-eps}
			\epsilon:=\frac{\mbox{Leb}(C)}{\sqrt{\eta}\,\sigma}\inf_{(x,y)\in C^2}\phi\left(\frac{y-x-\eta g(x)}{\sqrt{\eta}\sigma}\right),
		\end{align}
		$\nu(\cdot):=\mbox{Leb}(\cdot\,\cap C)/\mbox{Leb}(C)$, and~$\phi(x):=\frac{1}{\sqrt{2\pi}}\exp\{-x^2/2\}$ is the probability density function of standard normal distribution.
		\item[(2)] For any $\eta\in(0,1)$, the~Markov kernel $P_{\eta}$ is irreducible and strongly aperiodic.
		\item[(3)] Under Hypothesis \ref{hyp2}, the~state space $\mathbb{R}$ of the Markov chain $(\theta_k)_{k\geq0}$ is an $1$-small set.
	\end{enumerate}
\end{proposition}
\begin{pf}
	\begin{enumerate}
	\item[(1)] 	Suppose that $C$ is a compact subset of $\mathbb{R}$ such that $\mbox{Leb}(C)>0$. Recall that $P_{\eta}$ admits a Markov kernel $p_{\eta}(x,y)$, given~by
	$$p_{\eta}(x,y)=\frac{1}{\sqrt{\eta}\,\sigma}\phi\left(\frac{y-x-\eta g(x)}{\sqrt{\eta}\sigma} \right).$$
	Then for all $x\in C$ and $A\in\mathcal{B}(\R)$,		
	\begin{align*}
			P_{\eta}(x,A)=\int_A p_{\eta}(x,y) dy\geq\int_{A\cap C}p_{\eta}(x,y) dy\geq\,\epsilon\, \frac{\mbox{Leb}(A\cap C)}{\mbox{Leb}(C)},
	\end{align*}
	where $\epsilon\in(0,\,1]$ is defined in (\ref{Def-eps}). Evidently, $C$ is an $(1,\epsilon\nu)$-small set, with $\nu(\cdot)=\mbox{Leb}(\cdot\,\cap C)/\mbox{Leb}(C)$. Furthermore, $C$ is accessible since for all $x\in\R$,
	$$P_{\eta}(x,C)=\frac{1}{\sqrt{\eta}\,\sigma}\int_C\phi\left(\frac{y-x-\eta g(x)}{\sqrt{\eta}\sigma}\right)dy>0.$$
	\item[(2)] We have proved in (1) that any compact subset $C$ of $\R$ satisfying $\mbox{Leb}(C)>0$ is an accessible $(1,\mu)$-small set with $\mu(C)=\epsilon\nu(C)>0$. Hence, $P_{\eta}$ is irreducible and strongly aperiodic.
	\item[(3)] Assume Hypothesis \ref{hyp2}. Then we have that for a given non-empty set $A\in\mathcal{B}(\R)$, both $\inf_{x\in\R}(x+\eta g(x))$ and $\sup_{x\in\R}(x+\eta g(x))$ exist in $\R$; for a given $y\in A$,	
		$$
		\inf_{x\in\R}p_{\eta}(x,y)\geq\frac{1}{\sqrt{2\pi\eta}\sigma}\exp\left\{-\frac{1}{2\eta\sigma^2}\left[\left(y-\inf_{x\in\R}(x+\eta g(x))\right)^2\vee\left(y-\sup_{x\in\R}(x+\eta g(x))\right)^2\right]\right\}>0.
		$$
      Let $f(y)$ denote the function in the middle of the two inequalities above. Then for all $x\in\R$ and $A\in\mathcal{B}(\R)$,
      $$P_{\eta}(x,\,A)=\int_A p_{\eta}(x,y) dy\geq\mu(A),$$
      where $\mu(A)=\int_A f(y) dy$ is a Borel measure defined on $\mathcal{B}(\R)$. Therefore, the~state space $\R$ is an $(1,\mu)$-small set. \hfill$\square$
\end{enumerate}
\end{pf}

If $C$ is a subset of $\R$, let $\tau_C$ and $\sigma_C$ denote respectively the first hitting time and the first return time of $C$ for the chain $(\theta_n)_{n\geq 0}$, i.e.,
$$\tau_C=\inf\{n\geq0\,\big|\,\theta_n\in C\},$$
$$\sigma_C=\inf\{n\geq1\,\big|\, \theta_n\in C\}.$$
For a given $\eta\in(0,1)$, set
\begin{equation}\label{Def-b}
	b_{\eta}:=2\left(g^2(0)-L^2\right)\eta^2+\left(\frac{1}{2}K_1+\sigma^2+2c\right)\eta,
\end{equation}
where $L$, $K_1$ and $c$ are given from (\ref{assump1-1}), (\ref{assump1-2}) and (\ref{assump1-3}) respectively. Consider a function $\lambda(\eta)$ defined by
$$\lambda(\eta):=1-\frac{1}{2}K_1\eta+2L^2\eta^2,\quad \eta\in(0,1).$$
By Lemma \ref{lemApp} (1), there exists a constant $\eta_0\in(0,1)$ depending only on $K_1$ and $L$, such that for any $\eta\in(0,\eta_0]$, $0<\lambda(\eta)<1$. So by (\ref{Def-b}), we have for $\eta\in(0,\eta_0]$,
\begin{equation}\label{b_eta}
	b_{\eta}=1-\lambda(\eta)+2g^2(0)\eta^2+\sigma^2\eta+2c\eta>0.
\end{equation}
For $\eta\in(0,\eta_0]$, consider the sets defined as follows
$$D_{\eta}:=\left\{x\in\R\,\Big|\,|x|\leq\frac{2b_{\eta}}{K_1\eta^2}\right\}$$
and 
$$B_{\eta}:=\left\{x\in\R\,\Big|\,|x|\leq\frac{2b_{\eta}}{K_1\eta}\right\}.$$
Define a function $V(x)$ by  
$$V(x):=1+x^2,\quad x\in\R.$$

We have the following lemma. The~inequality (\ref{uniform-bound}) below shows that for any initial state $x\in \R$, the~quantity $\E_x\big(\beta_\eta^{\sigma_{D_{\eta}}}\big)$ is bounded from above uniformly in $\eta$.
\begin{lemma}\label{lem}
	Under Hypothesis \ref{hyp1}, the~following results~hold.
	\begin{enumerate}
		\item[(1)] For any $\eta\in(0,1)$, $D_{\eta}$ is an accessible $(1,\mu)$-small set with $\mu(D_{\eta})>0$.
		\item[(2)] There exists a constant $\eta_0\in(0,1)$ depending only on $K_1$ and $L$, such that for any $\eta\in(0,\eta_0]$ and $x\in\R$,
		$$\E_x\left(\beta_{\eta}^{\sigma_{D_\eta}}\right)\leq V(x)+b_{\eta}\beta_{\eta},$$
		where  $\beta_{\eta}>1$ is a constant depending on $\eta$. 
		
		Furthermore, for~any $\eta\in(0,\eta_0]$ and $x\in\R$,
		\begin{equation}\label{uniform-bound}
			\E_x\left(\beta_\eta^{\sigma_{D_{\eta}}}\right)\leq\E_x\left(\beta_{0}^{\sigma_{D_{0}}}\right)\leq  V(x)+b_{0}\beta_{0},
		\end{equation}
		where $D_0:=D_{\eta_0}$, $\beta_0:=\beta_{\eta_0}$ and $b_0:=b_{\eta_0}$.
	\end{enumerate}
\end{lemma}
\begin{pf}
	\begin{enumerate}
		\item[(1)] Since for any $\eta\in(0,1)$, $D_{\eta}$ is a compact subset of $\R$ and $\mbox{Leb}(D_{\eta})=\frac{4b_{\eta}}{K_1\eta^2}>0$, from~Proposition \ref{aperiodic} (1), we have that for any $\eta\in(0,1)$, $D_{\eta}$ is an accessible $(1,\mu)$-small set, with~$\mu=\epsilon\nu$, $\epsilon\in(0,\,1]$ and $\nu(D_{\eta})=1>0$, so that $\mu(D_{\eta})>0$.
		\item[(2)] Using (\ref{assump1-3}), (\ref{recursive}) and the first inequality of (\ref{eq-g-square}), we can prove (see also (A.2) in~\cite{LTX22}) that for any $\eta\in(0,1)$ and $x\in\R$, the~following inequality holds
		\begin{equation}\label{Drift-ine}
			P_{\eta}V(x)\leq\lambda(\eta) V(x)+b_{\eta}\mathbbm{1}_{B_{\eta}}(x).
		\end{equation}
		By (\ref{b_eta}), we have $\frac{b_\eta}{K_1\eta}<\frac{b_\eta}{K_1\eta^2}$, for~$\eta\in(0,\eta_0]$. Consequently, $B_\eta\subset D_\eta$ and $\mathbbm{1}_{B_\eta}\leq\mathbbm{1}_{D_\eta}$. From~(\ref{Drift-ine}), we obtain
		$$P_{\eta}V(x)\leq\lambda(\eta) V(x)+b_{\eta}\mathbbm{1}_{D_{\eta}}(x).$$
		Then according to {Proposition 4.3.3} (\rmnum{2}) in~\cite{DMPS18}, we have for any  $\eta\in(0,\eta_0]$ and $x\in \R$,
		\begin{align}\label{ine0}
			\E_x\left(\beta_{\eta}^{\sigma_{D_\eta}}\right)\leq V(x)+b_{\eta}\beta_{\eta},
		\end{align}
		where $\beta_{\eta}=1/\lambda(\eta)>1$.
		
		Applying Lemma \ref{lemApp} (1) again, since $\lambda(\eta)$ is decreasing in $\eta$, we have for any $\eta\in(0,\eta_0]$, $\beta_\eta\leq \beta_0$; from Lemma \ref{lemApp} (2), $D_{0}\subset D_{\eta}$, so \mbox{$\displaystyle\max_{\eta\in(0,\eta_0]}\sigma_{D_\eta}\leq \sigma_{D_0}$}. Therefore, by~(\ref{ine0}), we come to the result (\ref{uniform-bound}). \hfill$\square$
	\end{enumerate}
\end{pf}

Using the lemma above, we can obtain the following theorem, which shows the convergence of the kernel $P_{\eta}$ to its invariant probability measure $\pi_\eta$ under TV distance, for~any initial probability measure $\xi$ on $(\R,\,\mathcal{B}(\R))$.
\begin{theorem}\label{Thm1}
	Under Hypothesis $\ref{hyp1}$, there exists a constant $\eta_{0}\in(0,1)$ depending only on $K_1$ and $L$, such that for any $\eta\in (0,\,\eta_{0}]$, the~Markov kernel $P_{\eta}$ has a unique invariant probability measure $\pi_{\eta}$. Furthermore, there exist constants $\delta>1$, $\tau<\infty$,  $\beta_{0}>1$ and an accessible $(1,\mu)$-small set  $D_{0}$ (all depending only on $K_1$ and $L$) with $\mu(D_{0})>0$, such that
	\begin{equation}\label{geometric}
		\sup_{x\in D_0}\E_x\left(\beta_{0}^{\sigma_{D_{0}}}\right)<\infty,
	\end{equation}
	and for any initial probability measure $\xi$ on $(\R,\,\mathcal{B}(\R))$,
	\begin{equation}\label{series}
		\sum_{n=1}^{+\infty}\delta^n d_{TV}(\xi P_{\eta}^n,\pi_{\eta})\leq\tau \E_{\xi}\left(\beta_0^{\sigma_{D_0}}\right).
	\end{equation}
\end{theorem}
\begin{pf}
	From Lemma \ref{lem}, we have that $D_0$ is an accessible $(1,\mu)$-small set with $\mu(D_0)>0$ and there exists $\beta_0>1$ such that (\ref{geometric}) holds. Applying {Theorem 11.4.2} from~\cite{DMPS18} to the accessible $(1,\varepsilon\nu)$-small set $D_{0}$, we can obtain immediately the existence and uniqueness of the invariant probability measure $\pi_\eta$ for the Markov kernel $P_{\eta}$ and there exist constants $\delta>1$ and $\tau<\infty$ (both depending only on $K_1$ and $L$) such that for any initial probability measure $\xi$ on $(\R,\,\mathcal{B}(\R))$, (\ref{series}) holds. \hfill$\square$
\end{pf}

\begin{example}[AR(1) process]
	For the existence and uniqueness of invariant probability measure in Theorem \ref{Thm1}, one can take a simple example as follows. Let $g(x):=\left(1-\frac{1}{\eta}\right) x$, for~$x\in\R$; and $\sigma=1$. We can see that Hypothesis \ref{hyp1} is satisfied in this case. The~model $(\ref{recursive})$ becomes
	\begin{equation}\label{example}
		\theta_{k+1}=\eta\theta_k+\sqrt{\eta}\varepsilon_{k+1},\quad k\geq0,
	\end{equation}
	which is an $AR(1)$ model. Iterating (\ref{example}), we have
	$$\theta_{k}=\eta^k\theta_0+A_k,\quad k\geq1,$$
	where $A_k:=\sqrt{\eta}\sum_{i=0}^{k-1}\eta^i\varepsilon_{k-i}$. Let $\varepsilon_0$ be a standard normal random variable independent with $(\varepsilon_k)_{k\geq1}$. Consider the random variable $B_k:=\sqrt{\eta}\sum_{i=0}^{k-1}\eta^i\varepsilon_i$, for~any $k\geq1$. Since $0<\eta<1$, using martingale convergence theorem, we can obtain $\lim_{k\rightarrow+\infty}B_k=B_{\infty}$ a.s., where $B_\infty:=\sqrt{\eta}\sum_{i=0}^{+\infty}\eta^i\varepsilon_i$. Let $\pi_\eta$ be the probability distribution of $B_{\infty}$. Then we have that $\pi_\eta$ is the unique invariant probability measure for $P_\eta$.
\end{example}

\begin{remark}
	Theorem \ref{Thm1} improves the result in Lemma 2.3 of~\cite{LTX22} in the following two aspects. It shows not only the existence of invariant probability measure $\pi_{\eta}$, but~also the uniqueness under the same conditions, by~introducing the accessible small set $D_0$. It establishes furthermore the geometric rate of convergence of $\xi P_{\eta}^n$ towards $\pi_{\eta}$ under TV distance for any initial probability measure $\xi$ on $(\R,\mathcal{B}(\R))$, as~$n\rightarrow\infty$. See also the corollary below.
\end{remark}

\begin{corollary}\label{cor}
	Under Hypothesis $\ref{hyp1}$, there exist constants $\delta>1$ and $\eta_{0}\in(0,1)$ depending only on $K_1$ and $L$, such that for any initial probability measure $\xi$ on $(\R,\mathcal{B}(\R))$ and $\eta\in (0,\,\eta_{0}]$,
	\begin{equation*}
		d_{TV}\left(\xi P_{\eta}^n,\pi_\eta\right)=o\left(\delta^{-n}\right),\quad\text{as }n\rightarrow\infty.
	\end{equation*}
\end{corollary}
\section{{Geometric Rate of Convergence and Uniform Geometric Ergodicity of the Kernel $P_{ \eta}$}}\label{sec3}
Let $P$ be a Markov kernel on $X\times\mathcal{X}$. The~Markov kernel P is said to satisfy the geometric drift condition (see also {Definition 14.1.5}, \cite{DMPS18}) $D_{g}(V,\lambda, b)$, if~$V: X \rightarrow [1, \infty)$ is a measurable function, $\lambda\in[0,1)$, $b\in[0,\infty)$ and
\begin{equation}\label{drift}
	PV(x)\leq\lambda V(x)+b\mathbbm{1}_X(x).
\end{equation}

We obtain the following geometric rate of convergence for the Markov kernel $P_{\eta}$ to its invariant probability measure $\pi_\eta$, for~any initial state $x\in\R$.
\begin{theorem}\label{Thm2}
	Under Hypotheses $\ref{hyp1}$ and $\ref{hyp2}$, for~any $\eta\in (0,1)$, there exist $\delta\in(1,\infty]$ and $\tau<\infty$ such that for any $n\in\N$,
	\begin{equation}\label{Def_uniformlyGE}
		\sup_{x\in \R}\left\|P_{\eta}^n(x,\cdot)-\pi_{\eta}\right\|_{TV}\leq\tau\delta^{-n}
	\end{equation}
	and
	\begin{equation}\label{uniformlyGE}
		\sup_{x\in \R}d_{TV}\left(P_{\eta}^n(x,\cdot),\pi_{\eta}\right)\leq\left(\tau/2\right)\delta^{-n}.
	\end{equation}
\end{theorem}

 To prove the theorem above, we will need the following proposition, which gives equivalent conditions for uniform geometric ergodicity of $P_\eta$.
\begin{proposition}\label{uniformGE}
	For any $\eta\in(0,1)$, the~following statements are~equivalent.
	\begin{enumerate}
		\item [(1)] $P_{\eta}$ is uniformly geometrically ergodic, i.e.,~ $P_\eta$ admits an invariant probability measure $\pi_\eta$ such that there exist $\delta>1$ and $\tau<\infty$ satisfying for all $n\in\N$,
		\begin{equation}\label{ine-uniformGE}
			\sup_{x\in \R}\left\|P_{\eta}^n(x,\cdot)-\pi_{\eta}\right\|_{TV}\leq\tau\delta^{-n}.
		\end{equation}
		\item [(2)] $P_{\eta}$ is positive, aperiodic, and~there exist a small set $C$ and a constant $\delta>1$ such that
		\begin{equation*}
			\sup_{x\in \R}\E_x\left(\delta^{\sigma_C}\right)<\infty.
		\end{equation*}
		\item [(3)] The state space $\R$ is small.
	\end{enumerate}
\end{proposition}
\begin{pf} 
	(1) $\Longrightarrow$ (2). Suppose that $P_\eta$ is uniformly geometrically ergodic. We prove firstly that $P_\eta$ is irreducible; then as $P_{\eta}$ has an invariant probability measure $\pi_\eta$, $P_{\eta}$ is thus positive. Moreover, by~(\ref{ine-uniformGE}), there exist $\beta>1$ and $M<\infty$ such that for any $n\in\N$, $x\in \R$ and $A\in\mathcal{B}(\R)$, we have
	\begin{equation}\label{eq1}
		|P_{\eta}^n(x,\, A)-\pi_\eta(A)|\leq\|P^n_\eta(x,\cdot)-\pi_\eta\|_{TV}\leq M\beta^{-n}.
	\end{equation}
	Therefore, for~any $n\in\N$, $A\in\mathcal{B}(\R)$ and $x\in \R$,
	\begin{equation}\label{eq-origin}
		P^n_\eta(x,A)\geq\pi_\eta(A)-M\beta^{-n}.
	\end{equation}
	If $\pi_\eta(A)>0$, we can choose $n$ to be large enough such that $P^n_\eta(x,\,A)>0$, which implies that $P_\eta$ is~irreducible.
	
	Let $C$ be an accessible small set. For~$d>0$, define set $B_d=\{x\in\R\,|\,1(x)\leq d\}$. Since $\pi_\eta$ is an invariant probability measure, we get $\pi_\eta(C)>0$ and for all $x\in B_d$ we can choose $n$ to be large enough such that
	\begin{equation*}
		P_\eta^n(x,\,C)\geq\pi_\eta(C)-M\beta^{-n}\geq\pi_\eta(C)/2.
	\end{equation*}
	Therefore, for~all $d>0$, $B_d$ is also a small set by Lemma \ref{Lem9.1.7app1}.
	
	Since $\R=\{x\in\R\,|\,1(x)< \infty\}$ and $\pi_\eta(\R)=1$, we may choose $d_0$ to be large enough so that $\pi_\eta(B_d)>0$ for all $d\geq d_0$. Since $\pi_\eta$ is an invariant probability measure of $P_\eta$, the~set $B_d$ is accessible, for~any $d\geq d_0$. Applying (\ref{eq-origin}) for $A=B_d$, we may find $n$ to be large enough such that for any $m>n$, 
	$$\inf_{x\in B_d}P_\eta^m(x,B_d)\geq\pi_{\eta}(B_d)/2>0.$$ 
	This implies that the period of $P_\eta$ is $1$ and $P_\eta$ is thus~aperiodic. 
	
	On the one hand, letting $A=\R$ in (\ref{eq1}), one gets for any $k\in\N$ and $x\in \R$,
	\begin{equation}\label{eq3}
		P_\eta^k\,1(x)\leq M\beta^{-k}+1.
	\end{equation}
	We can thus choose $m\in\N^*$ to be large enough such that
	\begin{equation*}
		M\beta^{-m}\leq\lambda<1.
	\end{equation*}
	As a result, $P^{m}_\eta$ satisfies the geometric drift condition $D_g(1,\lambda,1)$, i.e.,~ $P^{m}_\eta\leq\lambda+1$. Now, set 
	$$V_0(x):=1+\lambda^{-1/m}P_\eta1(x)+\cdots+\lambda^{-(m-1)/m}P_{\eta}^{m-1}1(x),\quad x\in \R.$$
	Let $b:=\lambda^{-(m-1)/m}$. Then we have	
		\begin{align}
			P_\eta V_0(x)=&P_\eta1(x)+\lambda^{-1/m}P_\eta^21(x)+\cdots+\lambda^{-(m-1)/m}P_\eta^m1(x)\notag\\
			\leq& P_{\eta}1(x)+\lambda^{-1/m}P_\eta^21(x)+\cdots+\lambda^{-(m-2)/m}P_{\eta}^{m-1}1(x)+\lambda^{-(m-1)/m}(\lambda+1) \label{eq2}\\
			=&\lambda^{1/m}(\lambda^{-1/m}P_\eta1(x)+\lambda^{-2/m}P_\eta^21(x)+\cdots+\lambda^{-(m-1)/m}P_\eta^{m-1}1(x)+1)+b\notag\\
			=&\lambda^{1/m}V_0(x)+b.\notag
		\end{align}
	
	On the other hand, by~(\ref{eq3}), we have for any $x\in \R$,\vspace{-12pt}
      \begin{align*}
			1< V_0(x)\leq\left(M\sum_{k=1}^{m-1}\lambda^{-k/m}\beta^{-k}\right)+\sum_{k=0}^{m-1}\lambda^{-k/m}\leq(M+1)\sum_{k=0}^{m-1}\lambda^{-k/m}=(M+1)\frac{\lambda^{-1}-1}{\lambda^{-1/m}-1}.
		\end{align*}
	$P_\eta$ thus satisfies the geometric drift condition $D_g(V_0,\lambda^{1/m},\lambda^{-(m-1)/m})$, for~some $m\in\N^*$, $\lambda\in(0,1)$ and $V_0:\R\rightarrow[1,\infty)$ a measurable and bounded function on $\R$.
	
	Let $f(x):=(\tilde{\lambda}-\lambda^{1/m})V_0(x)$, for~$\tilde{\lambda}\in(\lambda^{1/m},\,1)$. From~(\ref{eq2}), for~any $\tilde{\lambda}\in(\lambda^{1/m},\,1)$ and $x\in \R$,
	\begin{equation}\label{eq-star}
		P_\eta V_0(x)+f(x)\leq\tilde{\lambda}V_0(x)+b.
	\end{equation} 
	Moreover, from~the above, we have for all $d>0$, the~set $B_d$ is small; and for all $d\geq d_0$, the~set $B_d$ becomes accessible. Therefore, these are also true for the set $\{x\in\R\,|\,V_0(x)\leq d\}$. Choose $\lambda_1\in(\tilde{\lambda},\, 1)$. Consider a set $C:=\{x\in\R\,|\,V_0(x)\leq d_1\}$, where $d_1\geq d_0\vee b\,(\lambda_1-\tilde{\lambda})^{-1}$; and hence $C$ is accessible and small.
	For $x\in C$, (\ref{eq-star}) yields $P_\eta V_0(x)+f(x)< \lambda_1V_0(x)+b$. While, for~$x\in C^c$, $-(\lambda_1-\tilde{\lambda})V_0(x)<-b$ and (\ref{eq-star}) imply
	\vspace{-4pt}
	$$P_\eta V_0(x)+f(x)<\lambda_1V_0(x)+b-(\lambda_1-\tilde{\lambda})V_0(x)<\lambda_1V_0(x).$$
	Equivalently, for~any $x\in\R$,
	\begin{equation}\label{eq5}
		P_\eta V_0(x)+f(x)\leq\lambda_1V_0(x)+b\mathbbm{1}_C.
	\end{equation}
	From the last inequality, with~measurable functions $V_0:\R\rightarrow[1,\infty]$ and $f: \R\rightarrow[1,\infty)$, some given $\delta=\lambda_1^{-1}>1$ and set $C$, we get for any $x\in C^c$,
	$$P_\eta V_0(x)+f(x)\leq\delta^{-1}V_0(x).$$
	According to ({14.1.4}) in {Proposition 14.1.2} \cite{DMPS18}, for~any $x\in\R$,
	$$\lambda_1^{-1}(\tilde{\lambda}-\lambda^{1/m})\,\E_x\left(\sum_{k=0}^{\sigma_C-1}\tilde{\lambda}^{-k}V_0(X_k)\right)\leq\lambda_1^{-1}\left[P_\eta V_0(x)+f(x)\right]\mathbbm{1}_C(x)+V_0(x)\mathbbm{1}_{C^c}(x),$$
	where $0\times\infty=0$ by convention on the right-hand side of the inequality. Using (\ref{eq5}), we~obtain
	\begin{align*}
		\lambda_1^{-1}(\tilde{\lambda}-\lambda^{1/m})\,\E_x\left(\sum_{k=0}^{\sigma_C-1}\delta^kV_0(X_k)\right)\leq&(\sup_CV_0+b\lambda_1^{-1})\mathbbm{1}_C(x)+V_0(x)\mathbbm{1}_{C^c}(x)\\
		\leq&\left(d_1+b\lambda_1^{-1}+1\right)V_0(x)=(d_1+b\delta+1)V_0(x).
	\end{align*}
	Therefore, for~any $x\in \R$,
	\begin{equation}\label{eq6}
		\E_x\left(\sum_{k=0}^{\sigma_C-1}\delta^kV_0(X_k)\right)\leq\tau V_0(x),
	\end{equation}
	where $\tau:=[(d_1+1)\delta^{-1}+b](\tilde{\lambda}-\lambda^{1/m})^{-1}<\infty$. Moreover, for~any $x\in \R$,
	\begin{align}\label{eq7}
		\E_x\left(\sum_{k=0}^{\sigma_C-1}\delta^kV_0(X_k)\right)=V_0(x)+\delta^{-1}\E_x\left(\sum_{i=2}^{\sigma_C}\delta^iV_0(X_{i-1})\right)\geq1+\delta^{-1}\E_x\left(\delta^{\sigma_C}\right).
	\end{align}
	Combining (\ref{eq6}) and (\ref{eq7}), we have that for any $x\in \R$, there exists $\delta\in(1,\infty)$ such that
	$$\E_x\left(\delta^{\sigma_C}\right)\leq\delta\left[\tau V_0(x)-1\right]<\infty.$$
	
	(2) $\Longrightarrow$ (3). For~any $d>0$, consider the set $D_d$ defined by
	$$D_d:=\{x\in \R\,|\,\E_x\left(\delta^{\tau_C}\right)<d\}.$$
	We will firstly prove that $D_d$ is small, for~any $d>0$. Note that 
	$$D_d=(D_d\cap C)\cup\left(D_d\cap C^c\right),$$ 
	and the union of two small sets remains small for irreducible and aperiodic kernel $P_\eta$. So it suffices to prove that both $D_d\cap C$ and $D_d\cap C^c$ are small. Since $C$ is small and $(D_d\cap C)\subseteq C$, it is evident that $D_d\cap C$ is also small. Suppose $x\in D_d\cap C^c$. By~Markov's inequality, for~all $k\in\N^*$,
	\begin{align*}
		\P_x(\sigma_C\geq k+1)=\P_x(\tau_C\geq k+1)=\P_x(\delta^{\tau_C}\geq\delta^{k+1})\leq\delta^{-(k+1)}\,\E_x\left(\delta^{\tau_C}\right)\leq d\delta^{-(k+1)}.
	\end{align*}
	Thus, for~$k$ that is sufficiently large,
	$$\inf_{x\in D_d\cap C^c }\P_x(\sigma_C\leq k)\geq\frac{1}{2}>0.$$
	Therefore, the~set $C$ is uniformly accessible from $D_d\cap C^c$. Since $P_\eta$ is irreducible, according to Lemma \ref{lem-uniformacess-app1} (1), we have that $D_d\cap C^c$ is~small. 
	
	Since
	$$\sup_{x\in\R}\E_x\left(\sigma^{\tau_C}\right)\leq\sup_{x\in\R}\E_x\left(\sigma^{\sigma_C}\right)<\infty,$$
	there exists $b>0$ such that $\R\subseteq D_b$, which is small from above. Therefore, the~state space $\R$ is~small.
	
	(3) $\Longrightarrow$ (1). Suppose that the state space $\R$ is small set. We will firstly show that $P_\eta$ is irreducible, recurrent and positive. Since $\R$ is small, there exists $m\in\N$ and nonzero measure $\mu$ on $(\R,\mathcal{B}(\R))$ such that for all $x\in\R$ and $A\in\mathcal{B}(\R)$,
	\begin{align*}
		P_\eta^m(x,A)\geq\mu_\eta(A),
	\end{align*}
	which implies for all $A\in\mathcal{B}(\R)$ such that $\mu_\eta(A)>0$, one has $P_\eta^m(x,A)>0$. Consequently, $P_\eta$ is irreducible. Since $\R$ is an accessible small set, $\sigma_\R=1$ $\P_x$-a.s., for~all $x\in\R$. Applying Lemma \ref{lem-uniformacess-app1} (2), we have that $P_\eta$ is recurrent. According to {Theorem 11.2.5} in~\cite{DMPS18}, $P_\eta$ admits an invariant probability measure $\pi_\eta$ satisfying $\pi_\eta(C)<\infty$ for every accessible set $C$; since $\R$ is accessible, this implies $\pi_\eta(\R)<\infty$, showing that $P_\eta$ is~positive.
	
	Next, we will prove that $P_\eta$ is aperiodic by contradiction. Suppose that $P_\eta$ is an irreducible Markov kernel with period $d>1$. According to Lemma \ref{Thm9.3.6app1}, there exist $C_0$, $C_1$, $\cdots$, $C_{d-1}$ of mutually disjoint accessible sets such that  for $i=0,\ldots,d-1$, and~for any $x\in C_i$, $P_\eta\left(x,C_{i+1(\text{mod }d)}\right)=1$ and hence
	\begin{align}\label{ine6}
		P_\eta(x,C_{(i+1)+1(\text{mod }d)})=0,
	\end{align}
	where additions are in the modulo sense, and~$\bigcup_{i=0}^{d-1}C_i$ is absorbing. Therefore, there exists $i_0\in\{0,\cdots,d-1\}$ such~that
	$$P_\eta^m(x,C_{i_0})\geq\mu_\eta(C_{i_0})>0,$$
	which contradicts with (\ref{ine6}).
	
	Finally, we can obtain (\ref{ine-uniformGE}) by applying Proposition \ref{prop-yinna1} below, with~$C=\R$. \hfill$\square$
\end{pf}

\begin{remark}
	From the proof above, we can see that the uniform geometric ergodicity of $P_\eta$ is connected to geometric drift condition (\ref{drift}). This allows us to control the quantity in the statement (2) and thus obtain the result in the statement (3). 
\end{remark}
\begin{proposition}\label{prop-yinna1}
	For any $\eta\in(0,1)$, suppose that $P_\eta$ admits an $(1,2\varepsilon\nu)$-small set $C$ with $\nu(C)=1$.  Assume that for some $\delta>1$,
	\begin{align}\label{ine15.1.0}
		\sup_{x\in C}\E_x\left(\sum_{k=0}^{\sigma_C}\delta^k\right)<\infty.
	\end{align} 
	Then there exists a constant $\beta>1$ such that for all initiate distributions $\xi$ on $(\R,\mathcal{B}(\R))$,
	\begin{align}\label{ine4}
		\sum_{n=0}^{\infty}\beta^n\|\xi P_\eta^n-\pi_\eta\|_{TV}<\infty.
	\end{align}
\end{proposition}

The proof of the proposition above will be given in the next~section.

\begin{pot2}
	The inequality (\ref{Def_uniformlyGE}) is an immediate consequence of Proposition \ref{aperiodic} (3) and Proposition \ref{uniformGE} (3). And~the inequality (\ref{uniformlyGE}) is obtained from (\ref{Def_uniformlyGE}) and the definition of TV distance (\ref{def_TVD}). \hfill$\square$
\end{pot2}
\section{Splitting Construction and Split Markov Kernel $ \check{P} _ \eta$}\label{sec4}
Let $\check{X}:=\R\times\{0,1\}$ and $(\check{X},\check{\mathfrak{X}})$ be a measurable space. In~this section, based on the kernel $P_\eta$, we  will construct a new Markov kernel $\check{P}_\eta$ on the extended state space $(\check{X},\check{\mathfrak{X}})$. The~detailed method of splitting construction for more general case can be found in {Section~11.1} of~\cite{DMPS18}.

Without loss of generality, we assume that $P_\eta$ admits an $(1,2\varepsilon\nu)$-small set $C$ with $\varepsilon\in(0,1)$ and $\nu(C)=1$. Let $b_{\varepsilon}$ be the Bernoulli distribution with success probability $\varepsilon$, given~by

$$b_\varepsilon:=(1-\varepsilon)\delta_{\{0\}}+\varepsilon\delta_{\{1\}}.$$
For any bounded and measurable function $f$ on $(\check{X},\check{\mathfrak{X}})$, define a function $\bar{f}_{\varepsilon}$ on $\R$ by

$$\bar{f}_{\varepsilon}(x):=[\delta_x\otimes b_\varepsilon]f=(1-\varepsilon)f(x,0)+\varepsilon f(x,1).$$
From the definition of $\bar{f}_{\varepsilon}$ above, we have for any measure $\xi$ on $(\R, \mathcal{B}(\R))$,
$$\xi(\bar{f}_{\varepsilon})=[\xi\otimes b_\varepsilon](f).$$
For $\eta\in(0,1)$, consider the residual kernel $R_\eta$ defined for $x\in\R$ and $A\in\mathcal{B}(\R)$ by

$$R_\eta(x,A):=\begin{cases}
	\frac{P_\eta(x,A)-\varepsilon\nu(A)}{1-\varepsilon},&\quad\text{if }x\in C;\\
	P_\eta(x,A),&\quad\text{if }x\notin C.
\end{cases}$$
Now, for~$\eta\in(0,1)$, define the split Markov kernel $\check{P}_\eta$ on $(\check{X},\check{\mathfrak{X}})$ as follows. For~$(x,d)\in\check{X}$ and $\check{A}\in\check{\mathfrak{X}}$, set

$$\check{P}_\eta(x,d;\check{A}):=Q_\eta(x,d;\cdot)\otimes b_\varepsilon(\check{A}),$$
where $Q_\eta$ is the Markov kernel on $\check{X}\times\mathcal{B}(\R)$ defined for all $B\in\mathcal{B}(\R)$ by
\vspace{-3pt}
$$Q_\eta(x,d;B):=\mathbbm{1}_C(x)\left[\mathbbm{1}_{\{0\}}(d)R_\eta(x,B)+\mathbbm{1}_{\{1\}}(d)\nu(B)\right]+\mathbbm{1}_{C^c}(x)P_\eta(x,B).$$
Equivalently, for~all bounded and measurable function $g$ on $(\R,\mathcal{B}(\R))$, one has
\vspace{-3pt}
\begin{align*}
	Q_\eta g(x,0)=&\mathbbm{1}_C(x)R_\eta g(x)+\mathbbm{1}_{C^c}(x)P_\eta g(x)=\begin{cases}
		\frac{\int_\R g(y)p_\eta(x,y)d y-\varepsilon\nu(g)}{1-\varepsilon},&\text{if }x\in~C\\
		\int_{\R}g(y)p_\eta(x,y)d y,&\text{if }x\notin C
	\end{cases};\\
	Q_\eta g(x,1)=&\mathbbm{1}_C(x)\nu (g)+\mathbbm{1}_{C^c}(x)P_\eta g(x)=\begin{cases}
		\nu(g),&\text{if }x\in~C\\
		\int_{\R}g(y)p_\eta(x,y)d y,&\text{if }x\notin C
	\end{cases},
\end{align*}
where $p_\eta(x,y)=(2\pi\eta\sigma^2)^{-1/2}\exp{\left(-\frac{(y-x-\eta g(x))^2}{2\eta\sigma^2}\right)}$ is the kernel density of $P_\eta$. It follows that for any bounded and measurable function $f$ on $(\check{X},\check{\mathfrak{X}})$,

$$\check{P}_\eta f(x,d)=Q_\eta\bar{f}_\varepsilon(x,d).$$

Let us describe now the canonical chain associated with the kernel $\check{P}_\eta$ on $\check{X}\times\check{\mathfrak{X}}$. For~probability measure $\check{\mu}$ on $\check{\mathfrak{X}}$, denote by $\check{\P}_{\check{\mu}}$ the probability measure on the canonical space $\left(\check{X}^{\N},{\check{\mathfrak{X}}}^{\otimes\N}\right)$ such that the coordinate process, denoted by $\{(\theta_k,D_k)\}_{k\geq1}$ and called split chain, is a Markov chain with initial distribution $\check{\mu}$ and Markov kernel $\check{P}_\eta$.

From the splitting construction above, $(D_n)_{n\geq0}$ is a sequence of i.i.d Bernoulli random variables with success probability $\varepsilon$ that is independent of $(\theta_n)_{n\geq1}$. Denote by $\left(\mathfrak{F}_k^\theta\right)_{k\geq1}$ the natural filtration of the process $(\theta_k)_{k\geq1}$. An~important property of the split chain $\{(\theta_k,D_k)\}_{k\geq1}$ is that if $\theta_0$ and $D_0$ are independent, then $\{(\theta_k,\mathfrak{F}_k^\theta)\}_{k\geq1}$ is a Markov chain with kernel $P_{\eta}$.

For a given $\eta\in(0,1)$, suppose that $P_\eta$ admits an $(1,\varepsilon\nu)$-small set $C$.  We have the following lemmas: Lemmas \ref{lem11.1.1}--\ref{prop11.1.3}.

\begin{lemma}\label{lem11.1.1}
	For any non-negative measure $\xi$ on $(\R,\,\mathcal{B}(\R))$,
	$$[\xi\otimes b_{\varepsilon}]\check{P}_{\eta}^n=\xi P_{\eta}^n\otimes b_\varepsilon.$$
\end{lemma} 

\begin{lemma}\label{prop11.1.2}
	For any probability measures on $(\R,\,\mathcal{B}(\R))$, $\{(\theta_k,\mathfrak{F}_k^\theta)\}_{k\geq1}$ is under $\check{\P}_{\xi\otimes b_\varepsilon}$ a Markov chain on $\R\times\,\mathcal{B}(\R)$ with initial distribution $\xi$ and Markov kernel $P_\eta$.
\end{lemma}

\begin{lemma}\label{prop11.1.3}
	If $\check{\lambda}_\eta$ is a non-negative measure on $(\check{X},\check{\mathfrak{X}})$ and is $\check{P}_\eta$-invariant, then $\check{\lambda}_{\eta}=\check{\lambda}_{\eta,0}\otimes b_\varepsilon$, where $\check{\lambda}_{\eta,0}$ is a non-negative measure on $(\R,\,\mathcal{B}(\R))$, defined by
	$$\check{\lambda}_{\eta,0}(A):=\check{\lambda}_{\eta}(A\times\{0,1\}),\quad A\in \mathcal{B}(\R).$$
	In addition, $\check{\lambda}_{\eta,0}$ is $P_\eta$-invariant.
\end{lemma}

For a given $\eta\in(0,1)$, suppose that $P_\eta$ admits an $(1,2\varepsilon\nu)$-small set $C$ with $\nu(C)=1$. We can obtain the following properties: Lemmas \ref{prop11.1.4} and \ref{prop-yinna}.

\begin{lemma}\label{prop11.1.4}
	Set $\check{\alpha}:=C\times\{1\}$ and $\check{C}:=C\times\{0,1\}$. Then the following results are~true:
	\begin{enumerate}
		\item [(1)] The set $\check{\alpha}$ is an aperiodic atom for the kernel $\check{P_\eta}$.
		\item [(2)] The set $\check{C}$ is small for the kernel $\check{P_\eta}$.
		\item [(3)] If $C$ is accessible, then the atom $\check{\alpha}$ is accessible for $\check{P_\eta}$, and~hence $\check{P}_\eta$ is irreducible.
		\item [(4)] For all $k\geq1$, $\check{P}_\eta^k\left(\check{\alpha},\check{\alpha}\right)=\varepsilon\nu P_\eta^{k-1}(C)$.
		\item [(5)] If $C$ is Harris recurrent for $P_\eta$, then for any probability measure $\xi$ on $(\R,\,\mathcal{B}(\R))$ satisfying $\P_{\xi}(\sigma_C<\infty)=1$, $\check{\P}_{\xi\otimes\delta_d}(\sigma_{\check{\alpha}}<\infty)=1$ for all $d\in\{0,1\}$. Moreover, if~$P_\eta$ is Harris recurrent, then $\check{P}_\eta$ is Harris recurrent.
		\item [(6)] If $C$ is accessible and $P_\eta$ admits an invariant probability measure $\pi_\eta$, then $\check{\alpha}$ is positive for $\check{P}_\eta$.
	\end{enumerate}
\end{lemma}

\begin{lemma}\label{prop-yinna}
	Assume that for some $\delta>1$,
	\begin{align}\label{ine15.1.1}
		\sup_{x\in C}\E_x\left(\sum_{k=0}^{\sigma_C}\delta^k\right)<\infty.
	\end{align}
	
	Then:
	\begin{enumerate}
		\item [(1)] There exist constants $\gamma\in(1,\delta)$ and $\tau<\tau_1<\infty$ such that
		\begin{align}\label{ine15.1.2}
			\sup_{(x,d)\in \check{C}}\check{\E}_{(x,d)}\left(\sum_{k=0}^{\sigma_{\check{\alpha}}}\gamma^k\right)\leq \tau\sup_{x\in C}\E_x\left(\sum_{k=0}^{\sigma_{C}-1}\delta^k\right),
		\end{align}
		and for any non-negative measure $\xi$ on $(\R,\,\mathcal{B}(\R))$,
		\begin{align}\label{ine15.1.3}
			\check{\E}_{\xi\otimes b_\varepsilon}\left(\sum_{k=0}^{\sigma_{\check{\alpha}}}\gamma^k\right)\leq \tau_1\E_{\xi}\left(\sum_{k=0}^{\sigma_{C}-1}\delta^k\right);
		\end{align}
		\item[(2)] $\check{P}_\eta$ admits a unique invariant probability measure $\pi_\eta\otimes b_\varepsilon$, where $\pi_\eta$ is the unique invariant probability measure of $P_\eta$.
	\end{enumerate}
\end{lemma}

\begin{pf}
	\begin{enumerate}
		\item [(1)] The condition (\ref{ine15.1.1}) implies that
		$$
		\inf_{x\in C}\P_x(\sigma_C<\infty)=1,
		$$
		so the set $C$ is Harris recurrent (see also {Proposition 4.2.5} (\rmnum{2}), \cite{DMPS18}). By~Lemma \ref{prop11.1.4} (5), for all $(x,d)\in\check{C}$, $\check{\P}_{(x,d)}(\sigma_{\check{C}}<\infty)=1$ and $\check{\P}_{(x,d)}(\sigma_{\check{\alpha}}<\infty)=1$. So we have for all $(x,d)\in \check{C}$ and $\gamma\in(1,\delta)$,
		\vspace{-3pt}
		\begin{align}
			\check{\E}_{(x,d)}(\gamma^{\sigma_{\check{\alpha}}})=&\gamma\,\check{\E}_{(x,d)}\left(\gamma^{\sigma_{\check{\alpha}}-1}\right)\notag\\
			\leq&\gamma\,\check{\E}_{(x,d)}\left(\sum_{k=0}^{\sigma_{\check{\alpha}}-1}\gamma^{k}\right)\label{ine15.1.4},
		\end{align}
		which implies that for any $(x,d)\in \check{C}$ and $\gamma\in(1,\delta)$,
		\begin{align}\label{ine15.1.4*}
			\check{\E}_{(x,d)}\left(\sum_{k=0}^{\sigma_{\check{\alpha}}}\gamma^{k}\right)\leq(\gamma+1)\,\check{\E}_{(x,d)}\left(\sum_{k=0}^{\sigma_{\check{\alpha}}-1}\gamma^{k}\right).
		\end{align}	
		On the other hand, 
		~for~any $x\in C$, we have by Lemma \ref{prop11.1.2},
		\begin{align}\label{eq15.1.5}
			\check{\E}_{\delta_x\otimes b_\varepsilon}\left(\sum_{k=0}^{\sigma_{\check{C}}-1}\delta^k\right)=\E_{x}\left(\sum_{k=0}^{\sigma_{C}-1}\delta^k\right).
		\end{align}
		Note that for any positive random variable $Y$,
		$$\sup_{(x,d)\in\check{C}}\check{\E}_{(x,d)}(Y)\leq\tau_\varepsilon\sup_{x\in C}\check{\E}_{\delta_x\otimes b_\varepsilon}(Y),$$
		with $\tau_\varepsilon:=\varepsilon^{-1}\vee(1-\varepsilon)^{-1}$. Applying the above inequality to $Y=\sum_{k=0}^{\sigma_{\check{C}}-1}\delta^k$, and~combining (\ref{eq15.1.5}) and the condition (\ref{ine15.1.1}), we get
		\begin{align}
			\sup_{(x,d)\in\check{C}}\check{\E}_{(x,d)}\left(\sum_{k=0}^{\sigma_{\check{C}}-1}\delta^k\right)\leq\tau_\varepsilon\sup_{x\in C}\check{\E}_{\delta_x\otimes b_\varepsilon}\left(\sum_{k=0}^{\sigma_{\check{C}}-1}\delta^k\right)=&\tau_\varepsilon\sup_{x\in C}\E_x\left(\sum_{k=0}^{\sigma_C-1}\delta^k\right)\label{eq15.0}\\
			<&\infty\label{eq15.1.5*}.
		\end{align}
		Similar to (\ref{ine15.1.4}), we can obtain
		\begin{align}\label{eq15.1.5-1}
			\check{\E}_{(x,d)}(\delta^{\sigma_{\check{C}}})\leq\delta\,\check{\E}_{(x,d)}\left(\sum_{k=0}^{\sigma_{\check{C}}-1}\delta^{k}\right).
		\end{align}
		Combining (\ref{eq15.1.5-1}) and (\ref{eq15.1.5*}), we obtain
		\begin{align*}
			\sup_{(x,d)\in\check{C}}\check{\E}_{(x,d)}\left(\delta^{\sigma_{\check{C}}}\right)<\infty.
		\end{align*}
		By Lemma \ref{prop11.1.4} (5), $\inf_{(x,d)\in\check{C}}\check{\P}_{(x,d)}(X_1\in\check{\alpha})>0$. Applying Theorem 14.2.3 in~\cite{DMPS18} with $A=\check{C}$, $B=\check{\alpha}$, $h=1$ and $q=1$, and~using (\ref{eq15.0}), we get there exist $\gamma\in(1,\delta)$, $\tau_0<\infty$ such that
		\begin{align*}
			\sup_{(x,d)\in\check{C}}\check{\E}_{(x,d)}\left(\sum_{k=0}^{\sigma_{\check{\alpha}}-1}\gamma^k\right)\leq&\tau_0\sup_{(x,d)\in\check{C}}\check{\E}_{(x,d)}\left(\sum_{k=0}^{\sigma_{\check{C}}-1}\delta^k\right)\\
			\leq&\tau_0\tau_\varepsilon\sup_{x\in C}\E_x\left(\sum_{k=0}^{\sigma_{C}-1}\delta^k\right).
		\end{align*}
		Combining the last inequality with (\ref{ine15.1.4*}), the~inequality (\ref{ine15.1.2}) is thus obtained with $\tau:=(\gamma+1)\,\tau_0\tau_\varepsilon<\infty$.
		
		Consider the shift operator $T$, defined by $T(\omega_0,\omega_1\ldots):=(\omega_1,\omega_2\ldots)$, for~any\linebreak   $\omega=(\omega_0,\omega_1\ldots)\in\R^{\N}$. Define inductively $T_0$ as the identity function, i.e.,~ $T_0(\omega)=\omega$, for~$\omega\in\R^{\N}$; and $T_n=T_{n-1}\circ T$, for~$n\geq1$.  Notice that $\sigma_{\check{\alpha}}\leq \sigma_{\check{C}}+\sigma_{\check{\alpha}}\circ T_{\sigma_{\check{C}}}$ on the event $\{\sigma_{\check{C}}<\infty\}$ and using Lemma \ref{prop11.1.2}, we get
		\begin{align*}
			\check{\E}_{\xi\otimes b_\varepsilon}\left(\sum_{k=0}^{\sigma_{\check{\alpha}}}\gamma^k\right)\leq&\check{\E}_{\xi\otimes b_\varepsilon}\left(\sum_{k=0}^{\sigma_{\check{C}}-1}\gamma^k\right)+\check{\E}_{\xi\otimes b_\varepsilon}\left(\sum_{k=\sigma_{\check{C}}}^{\sigma_{\check{\alpha}\circ T_{\sigma_{\check{C}}}}}\gamma^k\right)\\
			\leq&\check{\E}_{\xi\otimes b_\varepsilon}\left(\sum_{k=0}^{\sigma_{\check{C}}-1}\gamma^k\right)+\check{\E}_{\xi\otimes b_\varepsilon}\left(\gamma^{\sigma_{\check{C}}}\right)\sup_{(x,d)\in \check{C}}\check{\E}_{(x,d)}\left(\sum_{k=0}^{\sigma_{\check{\alpha}}}\gamma^k\right)\\
			\leq&\E_{\xi}\left(\sum_{k=0}^{\sigma_C-1}\gamma^k\right)\left[1+\gamma\sup_{(x,d)\in\check{C}}\check{\E}_{(x,d)}\left(\sum_{k=0}^{\sigma_{\check{\alpha}}}\gamma^k\right)\right].
		\end{align*}
		Therefore, by~(\ref{ine15.1.2}) and (\ref{ine15.1.1}), there exists $\tau_1\in (\tau,\infty)$ such that the inequality (\ref{ine15.1.3}) is satisfied.
		
		\item[(2)] Let $\eta\in(0,1)$ be fixed. From~Proposition \ref{aperiodic} (2), $P_\eta$ is irreducible and aperiodic. By~Lemma \ref{Coro9.2.14app1} and the condition (\ref{ine15.1.1}), the~small set $C$ is accessible. The~condition (\ref{ine15.1.1}) implies that $\sup_{x\in C}\E_x(\sigma_C)<\infty$. From~Lemma \ref{lem-uniformacess-app1} (3), we have $P_\eta$ is positive. Thus by Lemma \ref{prop11.1.4}, $\check{\alpha}$ is accessible, aperiodic, and~positive atom for the split kernel $\check{P}_\eta$. Using the inequality (\ref{ine15.1.2}), the~condition (\ref{ine15.1.1}) implies that there exists $\gamma\in(1,\delta)$ such that
		\begin{align}\label{ine1}
			\check{\E}_{\check{\alpha}}\left(\sum_{k=0}^{\sigma_{\check{\alpha}}}\gamma^k\right)<\infty.
		\end{align}
		{Applying} {Theorem 11.4.2} in~\cite{DMPS18} to $\check{P}_\eta$, we get it has a unique invariant probability measure $\check{\pi}_{\eta}$, which is expressed as $\pi_\eta\otimes b_\varepsilon$ by Lemma \ref{prop11.1.3}, where $\pi_\eta$ is the unique invariant probability measure for $P_\eta$. \hfill$\square$
	\end{enumerate}
\end{pf}

Using the lemmas above, we can prove Proposition \ref{prop-yinna1} as follows.
\begin{pot3}
	Lemma \ref{lem11.1.1} implies for $n\geq1$,
	$$\|\xi P_\eta^n-\pi_\eta\|_{TV}\leq\|(\xi\otimes b_\varepsilon)\check{P}_\eta^n-\pi_\eta\otimes b_\varepsilon\|_{TV}.$$
	So we have
	\begin{align}\label{ine5}
		\sum_{n=1}^\infty\beta^n\|\xi P_\eta^n-\pi_\eta\|_{TV}\leq\sum_{n=1}^\infty\beta^n\|(\xi\otimes b_\varepsilon)\check{P}_\eta^n-\pi_\eta\otimes b_\varepsilon\|_{TV}.
	\end{align}
	In Lemma \ref{prop-yinna}, it is proven that $\check{\alpha}$ is an accessible, aperiodic, and~positive atom for $\check{P}_\eta$, which admits unique invariant probability measure $\pi_\eta\otimes b_\varepsilon$; moreover, it is proven that (\ref{ine1}) is true. Now, applying Lemma \ref{Thm13.4.3app2} to $\check{P}_\eta$ with $\alpha=\check{\alpha}$, we obtain that there exist $\beta\in(1,\gamma)$ and $\tau<\infty$ such~that
	$$	\sum_{n=1}^\infty\beta^n\|(\xi\otimes b_\varepsilon)\check{P}_\eta^n-\pi_\eta\otimes b_\varepsilon\|_{TV}\leq\tau\check{\E}_{\xi\otimes b_\varepsilon}\left(\sum_{n=1}^{\sigma_{\check{\alpha}}}\gamma^n\right).$$
	Combining this inequality with (\ref{ine15.1.3}) and (\ref{ine5}) in Lemma \ref{prop-yinna} and the condition (\ref{ine15.1.1}), we can obtain the desired inequality (\ref{ine4}). \hfill$\square$
\end{pot3}
\section{Conclusions}\label{sec5}
This work investigates the convergence rate of EM scheme $(\theta_k)_{k\geq0}$ to the invariant measure under total variation distance for solution approximations of SDE given by (\ref{SDE}). To~this end, the~ergodic properties of the EM scheme, as~a non-atomic Markov chain, are studied. It turns out that under Hypothesis \ref{hyp1}, $(\theta_k)_{k\geq0}$ is irreducible, strongly aperiodic (see also Proposition \ref{aperiodic}), and~admits a unique invariant probability measure $\pi_\eta$ (see also Theorem \ref{Thm1}). If~in addition Hypothesis \ref{hyp2} is satisfied, $(\theta_k)_{k\geq0}$ is uniformly geometrically ergodic, and~converges to $\pi_\eta$ under TV distance at least in geometric rate, which is independent of the step size $\eta$ (see also Theorem \ref{Thm2}). To~show the uniform geometric ergodicity of the chain, we studied its equivalent conditions (see also Proposition \ref{uniformGE}) with the approach of splitting~construction.

Let us describe some connections and special highlights of the present work with classical ergodic theory for Markov chains on general state space. In~the book of Meyn and Tweedie (\cite{MT09}), the~theory on geometric ergodicity for $\psi$-irreducible Markov chain, with~Markov kernel $P$ and state space $X$, is introduced in Chapter 15. Our geometric ergodicity results in Theorem \ref{Thm1} are similar to those under $f$-norm in Theorem 15.4.1 of~\cite{MT09}, when $f=\mathbbm{1}_{X}$. Their results hold under the following ``Minorization Condition'' of the chain: for some $\delta>0$, some $C\in\mathcal{B}(X)$ and some probability measure $\nu$ with $\nu(C^c)=0$ and $\nu(C)=1$,  
\begin{equation*}
	P(x,A)\geq\delta\mathbbm{1}_C(x)\nu(A),\quad A\in\mathcal{B}(X),\; x\in X.
\end{equation*}

While, our geometric ergodicity results in Theorem \ref{Thm1}, obtained by using {Theorem 11.4.2} in~\cite{DMPS18}, only require the condition of existence of an $(1,\varepsilon\nu)$-small set $D$, with~$\varepsilon\in (0,1)$ and $\nu(D)=1$, which is ensured in our case because it is shown in Proposition 3.1 (3) that the state space is small. In~this article, we applied a new splitting construction approach introduced by Dedecker and Gou\"{e}zel~\cite{DG15}, which is different from the two classical splitting approaches in~\cite{MT09}: Nummelin splitting technique, due to~\cite{N84}, and~random renewal time approach, due to~\cite{AN78}. Thanks to this new splitting construction, it makes our proof of the equivalent conditions for the uniform geometric ergodicity of $P_\eta$ in Proposition \ref{uniformGE} differ from the proof of the existing similar results in {Theorems 16.0.2 and 16.2.2} of~\cite{MT09}. Moreover, due to this new approach, we have managed to establish an important link between $1$-geometrically recurrent small set as described in (\ref{ine15.1.0}) and the convergence of the series involving TV norm in (\ref{ine4}) for the Markov kernel $P_\eta$.

With the new concepts related to non-atomic Markov chains, our results provide further properties for an EM scheme with constant step size of SDEs. In~the future, EM scheme with non-constant step sizes and algorithm with high performance such as convergence rate in a controlled manner, can be considered for applications. Also, the~convergence behavior of the EM scheme in higher-dimensional spaces under different distances can be considered as a future research direction. 
\appendix
\section{Some Properties Related to Small Sets and Irreducible Markov Kernels}\label{sec5-1}
In this section, we will introduce some basic properties related to Markov kernels that admit accessible small sets. Readers may refer to Chapter 9 in the book of Douc~et~al.~\cite{DMPS18} for the~details.

Suppose that $(X_n)_{n\geq0}$ is a Markov chain living on the state space $X$, with~Markov kernel $P$ on a measurable space $(X,\mathfrak{X})$, where $\mathfrak{X}$ is a $\sigma$-algebra generated by $X$. In~the sequel, let $\P_x$ denote the probability measure on the canonical space $\left(X^{\N},\mathfrak{X}^{\otimes\N}\right)$ of the chain $(X_n)_{n\geq0}$, with~the initial state $X_0=x$. And~$\E_x$ denotes the expectation under the probability measure $\P_x$. If~$B\in\mathfrak{X}$, let $N_B$ denote the number of visits of $(X_n)_{n\geq0}$ to the set $B$, defined as $N_B=\sum_{k=0}^\infty\mathbbm{1}_B(X_k)$. Define respectively the first return time $\sigma_B$ of the set $B$ and the expected number $U(x,B)$ of visits to $B$ starting from $x$ by
$$\sigma_B=\inf\{n\geq1\,|\, X_n\in B\} $$
and
$$U(x,B)=\E_x[N_B]=\sum_{k=0}^{\infty}P^k(x,B).$$

\begin{definition}[Accessible set and uniform accessibility,~\cite{DMPS18}]
	\begin{enumerate}
		\item [(1)] \label{accessible} A set $A\in\mathfrak{X}$ is said to be   accessible  if for all $x\in X$, there exists an integer $n\geq1$ such that $P^n(x,A)>0$.
		\item[(2)]\label{def-uniformaccess} A set $B$ is   uniformly accessible from $A$, if~there exists $m\in\N^*$ such that
		$$\inf_{x\in A}\P_x(\sigma_B\leq m)>0.$$
	\end{enumerate}
\end{definition}
\begin{definition}[Full set,~\cite{DMPS18}]
	A set $F\in\mathfrak{X}$ is said to be  full if $F^c$ is not accessible.
\end{definition}
\begin{definition}[Atom and small set,~\cite{DMPS18}]\label{smallset}
	\begin{enumerate}
		\item[(1)] A set $\alpha\in\mathfrak{X}$ is called an  atom, if~there exists a probability measure $\nu$ on on $(X,\mathfrak{X})$ such that for all $x\in \alpha$ and $A\in \mathfrak{X}$,
		$$P(x,A)=\nu(A).$$
		\item[(2)] A set $C\in\mathfrak{X}$ is called a  small set if there exists positive integer $m$ and a nonzero measure $\mu$ on $(X,\mathfrak{X})$ such that for all $x\in C$ and $A\in \mathfrak{X}$,
		\begin{equation}\label{small}
			P^m(x,A)\geq\mu(A).
		\end{equation}
		Then the set $C$ is said to be an  $(m,\mu)$-small set.
	\end{enumerate}
\end{definition}
\begin{remark}\label{remark}
	\begin{enumerate}
		\item[(1)] From the definition above, it can be seen that atom is a particular small set satisfying the equality in (\ref{small}) with $m=1$ and $\mu(X)=1$, where $\mu$ becomes a probability measure.
		\item[(2)] The condition (\ref{small}) implies that $\mu$ is a finite measure satisfying $0<\mu(X)\leq1$. Hence it can be written $\mu=\epsilon\nu$ with $\epsilon=\mu(X)$, and~$\nu(\cdot)=\mu(\cdot)/\mu(X)$ is a probability measure on $(X,\mathfrak{X})$. If~$\epsilon=1$, then the equality in (\ref{small}) must hold, and~thus $A$ is an atom. 
	\end{enumerate}
\end{remark}
\begin{definition}[Irreducible kernel,~\cite{DMPS18}]\label{irreducible} A Markov kernel $P$ is said to be   irreducible if it admits an accessible small set.
\end{definition}
\begin{definition}[Recurrent and Harris recurrent,~\cite{DMPS18}]
	\begin{enumerate}
		\item [(1)] A set $A\in\mathfrak{X}$ is said to be  recurrent if $U(x,A) = \infty$ for all $x\in A$; it is said to be  Harris recurrent if $\P_x(N_A=\infty)=1$, for~all $x\in A$.
		\item[(2)] The kernel $P$ is said to be  recurrent if every accessible set is recurrent; it is said to be  Harris recurrent if every accessible set is Harris recurrent.
	\end{enumerate}
\end{definition}
\begin{definition}[Strongly aperiodic small set,~\cite{DMPS18}] An $(m,\mu)$-small set $C$ is said to be  strongly aperiodic, if~$m=1$ and $\mu(C)>0$.
\end{definition}
\begin{definition}[Period, aperiodicity and strong aperiodicity~\cite{DMPS18}]\label{saKernal}
	\begin{enumerate}
		\item [(1)] The common period of all accessible small sets is called the period of the kernel $P$.
		\item [(2)] If the period is equal to 1, the~kernel $P$ is said to be  aperiodic.
		\item [(3)] If there exists an accessible $(1,\mu)$-small set $C$ with $\mu(C)>0$, the~kernel $P$ is said to be  strongly aperiodic.
	\end{enumerate}
\end{definition}
\begin{definition}[Positive and null recurrent atom, positive and null Markov kernel,~\cite{DMPS18}]\label{Def_P} 
	\begin{enumerate}
		\item[(1)] An atom $\alpha$ is said to be  positive, if~$\E_{\alpha}(\sigma_\alpha)<\infty$; it is said to be  null recurrent, if~it is recurrent and $\E_{\alpha}(\sigma_\alpha)=\infty$.
		\item[(2)] If $P$ is irreducible and admits an invariant probability measure $\pi$, the~		Markov kernel $P$ is called  positive. If~$P$ does not admit such a measure, then
		$P$ is called  null.
	\end{enumerate}
\end{definition}
\begin{definition}[Uniformly geometrically ergodic,~\cite{DMPS18}]\label{Def_UGE}
	A Markov kernel $P$ on $X\times\mathfrak{X}$ is said to be  uniformly geometrically ergodic, if~it admits an invariant probability measure $\pi$ such that there exist $\delta>1$ and $\tau<\infty$ satisfying for any $n\in\N$ and any $x\in X$,
	$$\|P^n(x,\cdot)-\pi\|_{TV}\leq \tau\delta^{-n}.$$
\end{definition}
\begin{definition}[Invariant measure,~\cite{DMPS18}]
	A nonzero measure $\mu$ is said to be  invariant if it is $\sigma$-finite and $\mu P=\mu$.
\end{definition}

The lemmas below are useful to prove Proposition \ref{uniformGE}.

\begin{lemma}[{Lemma 9.1.7} (\rmnum{2}), \cite{DMPS18}]\label{Lem9.1.7app1}
	Let $C$ be an $(m,\mu)$-small set and $D\in\mathfrak{X}$. If~there exists $n\in\N$ such that $\inf_{x\in D}P^n(x, C)\geq\delta$, then $D$ is an $(n+m, \delta\mu)$-small set.
\end{lemma}

\begin{lemma}[{Corollary 9.2.14}, \cite{DMPS18}]\label{Coro9.2.14app1}
	Suppose that $P$ is irreducible. Let $r$ be a positive increasing sequence such that $\lim_{n\rightarrow\infty}r(n)=\infty$ and $A\in \mathfrak{X}$, $A\neq\emptyset$. Assume that $sup_{x\in A}\E_x[r(\sigma_A)] < \infty$. Then the set
	$$\{x\in X\;|\;\E_x[r(\sigma_A)]<\infty\}$$
	is full and absorbing, and~$A$ is accessible.
\end{lemma}

\begin{lemma}[{Theorem 9.3.6}, \cite{DMPS18}]\label{Thm9.3.6app1} Suppose that $P$ is an irreducible Markov kernel with period $d$. Then there exists a sequence $C_0,C_1,\ldots,C_{d-1}$ of mutually disjoint accessible sets such that for $i=0,\ldots,d-1$ and $x\in C_i$, $P\left(x,C_{i+1[mod\;d]}\right)=1$. Consequently, $\bigcup_{i=0}^{d-1}C_i$ is absorbing.
\end{lemma}

\begin{lemma}\label{lem-uniformacess-app1}
	Suppose that $P$ is irreducible and aperiodic. Then the following statements are~true. 
	\begin{enumerate}
		\item [(1)] Let $C$, $D\in\mathfrak{X}$. If~$D$ is small and uniformly accessible from $C$, then $C$ is also small.
		\item [(2)] $P$ is recurrent if and only if it admits an accessible recurrent small set.
		\item [(3)] If there exists a small set $C$, such that $\sup_{x\in C}\E_x(\sigma_C)<\infty$, then $P$ is positive.
	\end{enumerate}
\end{lemma}
The statements (1), (2) and (3) above can be proved by applying {Lemma 9.4.7} (\rmnum{2}) and {Theorem 9.4.10, Theorem 10.1.2 and Theorem 9.4.10, Corollary 11.2.9 and Theorem 9.4.10} respectively in~\cite{DMPS18}.
\section{Rate of Convergence in TV Norm for Atomic Markov~Chains}\label{sec5-2}
Suppose now that $\alpha$ is an atom for the kernel $P$. If~a function $h$ defined on $X$ is constant on $\alpha$, then we write $h(\alpha)$ instead of $h(x)$ for all $x\in\alpha$. With~this convention, for~every positive $\mathfrak{X}^{\N}$-measurable random variable $Y$ such that $\E_x(Y)$ is constant on $\alpha$, we write $\E_{\alpha}(Y)$ instead of $\E_{x}(Y)$, for~any $x\in\alpha$.

The following result shows the convergence rate of $P^n$ to its invariant probability measure $\pi$ under TV norm, if~the atom $\alpha$ is $1$-geometrically recurrent (see also {Definition 14.4.1}, \cite{DMPS18}).
\begin{lemma}\label{Thm13.4.3app2} Let $\alpha$ be an accessible, aperiodic, and~positive atom. Denote by $\pi$ the unique invariant probability measure. Assume that there exists $\gamma > 1$ such that
	$$
	\E_\alpha\left[\sum_{n=1}^{\sigma_\alpha}\gamma^n\right]<\infty.
	$$
Then there exist $\beta\in(1,\gamma)$ and a constant $\tau<\infty$ such that for every probability measure $\xi$ on $(X,\mathfrak{X})$,
	$$\sum_{n=1}^\infty\beta^n\|\xi P^n-\pi\|_{TV}\leq\tau\E_\xi\left[\sum_{n=1}^{\sigma_\alpha}\gamma^n\right].$$
\end{lemma}
The lemma above can be obtained by applying {Theorem 13.4.3} in~\cite{DMPS18} with $f=\mathbbm{1}_X$ therein.

\section{Some Auxiliary Result}\label{sec5-3}
Recall that
\begin{equation}\label{Def2-b}
	b_{\eta}:=2\left(g^2(0)-L^2\right)\eta^2+\left(\frac{1}{2}K_1+\sigma^2+2c\right)\eta,\quad \eta\in(0,1).
\end{equation}

In this part, we will conduct an analysis on the functions $\lambda(\eta)$ and $f_1(\eta)$ defined respectively as follows. For~any $\eta\in(0,1)$, $\lambda(\eta):=1-\frac{K_1}{2}\eta+2L^2\eta^2$ and
\begin{equation}\label{Def-f1}
	f_1(\eta):=\frac{2b_\eta}{K_1\eta^2}.
\end{equation}

We have the following result.
\begin{lemma}\label{lemApp}
	\begin{enumerate}
		\item [(1)] There exists a constant $\eta_{0}\in(0,1)$ depending only on $K_1$ and $L$ such that for all $\eta\in(0,\eta_0]$, the~function $\eta\mapsto\lambda(\eta)$ is decreasing and $0<\lambda(\eta)<1$.
		\item [(2)] The function $\eta\mapsto f_1(\eta)$ is decreasing on $(0,1)$.
	\end{enumerate}
\end{lemma}
\begin{pf}
	\begin{enumerate}
		\item [(1)]	The function $\lambda(\eta)$ is quadratic polynomial and can be written as
		$$\lambda(\eta)=2L^2\left(\eta-\frac{K_1}{8L^2}\right)^2+1-\frac{K_1^2}{32L^2}.$$
		Taking into account that $\lambda(0)=1$ and $\lambda(\eta)$ is symmetric w.r.t $\eta=K_1/8L^2>0$, we can immediately obtain the result.
		\item[(2)] By (\ref{Def2-b}) and (\ref{Def-f1}), we have for $\eta\in(0,1)$,
		$$f_1(\eta)=\frac{4\left(g^2(0)-L^2\right)}{K_1}+\frac{K_1+2\sigma^2+4c}{K_1\eta}$$
		and
		$$f'_1(\eta)=-\frac{K_1+2\sigma^2+4c}{K_1\eta^2}<0.$$
		Consequently, the~function $f_1(\eta)$ is decreasing on $(0,1)$. \hfill$\square$
	\end{enumerate}
\end{pf}






\end{document}